\documentclass[10pt]{article}
\usepackage{amssymb,latexsym}
\input epsf

\title{On the enumeration of positive cells in generalized cluster 
complexes and Catalan hyperplane arrangements}

\author{Christos~A.~Athanasiadis and Eleni~Tzanaki\\
Department of Mathematics\\[-0.6ex]
University of Crete\\[-0.6ex]
71409 Heraklion, Crete, Greece\\[-0.6ex]
\small \texttt{caath@math.uoa.gr, etzanaki@math.uoc.gr}}

\date{\small March 31, 2005. Revised; September 30, 2005. }

\newcommand{\Cat}{{\rm Cat}}

\newcommand{\NN}{{\mathbb N}}
\newcommand{\ZZ}{{\mathbb Z}}
\newcommand{\aA}{{\mathcal A}}

\newcommand{\fF}{{\mathcal F}}
\newcommand{\gG}{{\mathcal G}}
\newcommand{\iI}{{\mathcal I}}
\newcommand{\jJ}{{\mathcal J}}
\newcommand{\rR}{{\mathcal R}}

\newcommand{\ppP}{{\bf P}}
\newcommand{\qqQ}{{\bf Q}}
\renewcommand{\to}{\rightarrow}

\newcommand{\sm}{{\setminus}}
\newcommand{\qed}{$\hfill \Box$}

\begin{document}
\maketitle

\newtheorem{theorem}{Theorem}[section]
\newtheorem{proposition}[theorem]{Proposition}
\newtheorem{corollary}[theorem]{Corollary}
\newtheorem{definition}[theorem]{Definition}
\newtheorem{remark}[theorem]{Remark}
\newtheorem{lemma}[theorem]{Lemma}
\newtheorem{example}[theorem]{Example}
\newtheorem{examples}[theorem]{Examples}
\newtheorem{conjecture}[theorem]{Conjecture}
\newtheorem{question}[theorem]{Question}

\begin{abstract}
Let $\Phi$ be an irreducible crystallographic root system with Weyl 
group $W$ and coroot lattice $\check{Q}$, spanning a Euclidean space 
$V$. Let $m$ be a positive integer and $\aA^m_\Phi$ be the arrangement 
of hyperplanes in $V$ of the form $(\alpha, x) = k$ for $\alpha \in 
\Phi$ and $k = 0, 1,\dots,m$. It is known that the number $N^+ (\Phi, 
m)$ of bounded dominant regions of $\aA^m_\Phi$ is equal to the number 
of facets of the positive part $\Delta^m_+ (\Phi)$ of the generalized 
cluster complex associated to the pair $(\Phi, m)$ by S.~Fomin and 
N.~Reading. 

We define a statistic 
on the set of bounded dominant regions of $\aA^m_\Phi$ and conjecture 
that the corresponding refinement of $N^+ (\Phi, m)$ coincides with the 
$h$-vector of $\Delta^m_+ (\Phi)$. We compute these refined numbers for 
the classical root systems as well as for all root systems when $m=1$
and verify the conjecture when $\Phi$ has type $A$, $B$ or $C$ and when 
$m=1$. 
We give several combinatorial interpretations to these numbers in terms 
of chains of order ideals in the root poset of $\Phi$, orbits of the 
action of $W$ on the quotient $\check{Q} / \, (mh-1) \, \check{Q}$ and 
coroot lattice points inside a certain simplex, analogous to the ones 
given by the first author in the case of the set of all dominant regions 
of $\aA^m_\Phi$. We also provide a dual interpretation in terms of order 
filters in the root poset of $\Phi$ in the special case $m=1$.
\end{abstract}

\section{Introduction and results}
\label{intro}

Let $V$ be an $\ell$-dimensional Euclidean space, with inner product 
$( \ , \ )$. Let $\Phi$ be a (finite) irreducible crystallographic 
root system spanning $V$ and $m$ be a fixed nonnegative integer. We 
denote by $\aA_\Phi^m$ the collection of hyperplanes in $V$ defined 
by the affine equations $(\alpha, x) = k$ for $\alpha \in \Phi$ and 
$k = 0, 1,\dots,m$, known as the $m$th \textit{extended Catalan 
arrangement} associated to $\Phi$. Thus $\aA_\Phi^m$ is invariant 
under the action of the Weyl group $W$ associated to $\Phi$ and 
reduces to the Coxeter arrangement $\aA_{\Phi}$ for $m=0$. Let 
$\Delta^m (\Phi)$ denote the generalized cluster complex associated 
to the pair $(\Phi, m)$ by S.~Fomin and N.~Reading \cite{FR2}. This 
is a simplicial complex which reduces to the cluster complex $\Delta 
(\Phi)$ of S.~Fomin and A.~Zelevinsky \cite{FZ} when $m=1$. It 
contains a natural subcomplex, called the positive part of 
$\Delta^m (\Phi)$ and denoted by $\Delta^m_+ (\Phi)$, as an induced 
subcomplex. The complex $\Delta^m (\Phi)$ was also studied independently
by the second author \cite{Tz} when $\Phi$ is of type $A$ or $B$; see 
Section \ref{pre} for further information and references. 

The Weyl group $W$ acts on the coroot lattice $\check{Q}$ of $\Phi$ 
and its dilate $(mh-1) \, \check{Q}$, where $h$ denotes the Coxeter 
number of $\Phi$. Hence $W$ acts also on the quotient $T_m = \check{Q} 
/ \, (mh-1) \, \check{Q}$.
For a fixed choice of a positive system $\Phi^+ \subseteq \Phi$, 
consider the partial order on $\Phi^+$ defined by letting $\alpha \le 
\beta$ if $\beta - \alpha$ is a nonnegative linear combination of 
positive roots, known as the \emph{root poset} of $\Phi$. An 
\textit{order filter} or \textit{dual order ideal} in $\Phi^+$ is a 
subset $\iI$ of $\Phi^+$ such that $\alpha \in \iI$ and $\alpha \le 
\beta$ in $\Phi^+$ imply $\beta \in \iI$. The filter $\iI$ is called 
\emph{positive} if it does not contain any simple root.

The following theorem connects the objects just discussed. Parts (i), 
(ii) and (iii) appear in \cite[Corollary 1.3]{Ath1}, \cite[Proposition
2.13]{FR2} and \cite[Theorem 7.4.2]{Ha}, 
respectively. The last statement was found independently in \cite{Ath1, 
Pa, So}.

\begin{theorem} \mbox{\rm (\cite{Ath1, FR2, Ha})}
Let $\Phi$ be an irreducible crystallographic root system of rank $\ell$
with Weyl group $W$, Coxeter number $h$ and exponents $e_1, 
e_2,\dots,e_\ell$. Let $m$ be a positive integer and let
\[ N^+ (\Phi, m) = \prod_{i=1}^{\ell} \frac{e_i + mh - 1}{e_i + 1}. \]
The following are equal to $N^+ (\Phi, m)$:
\begin{enumerate}
\itemsep=0pt
\item[{\rm (i)}] the number of bounded regions of $\aA^m_\Phi$ which lie in the 
fundamental chamber of $\aA_\Phi$,
\item[{\rm (ii)}] the number of facets of $\Delta^m_+ (\Phi)$ and
\item[{\rm (iii)}] the number of orbits of the action of $W$ on $\check{Q} / \, 
(mh-1) \, \check{Q}$.
\end{enumerate} 
Moreover, for $m=1$ this number is equal to the number of positive 
filters in the root poset of $\Phi$.
\label{thm0}
\end{theorem}

The purpose of this paper is to define and study a refinement of the number
$N^+ (\Phi, m)$ and prove that it has similar properties with the one defined
by the first author \cite{Ath2} for the total number 
\[ N (\Phi, m) = \prod_{i=1}^{\ell} \frac{e_i + mh + 1}{e_i + 1} \]
of regions of $\aA^m_\Phi$ in the fundamental chamber of $\aA_\Phi$. To be 
more precise let $H_{\alpha, k}$ 
be the affine hyperplane in $V$ defined by the equation $(\alpha, x) = k$ 
and $A_\circ$ be the fundamental alcove of the affine Weyl arrangement 
corresponding to $\Phi$. A \emph{wall} of a region $R$ of $\aA^m_\Phi$ is 
a hyperplane in $V$ which supports a facet of $R$. For $0 \le i \le \ell$ 
we denote by $h_i (\Phi, m)$ the number of regions $R$ of $\aA_\Phi^m$ in 
the fundamental chamber of $\aA_\Phi$ for which exactly $\ell - i$ walls of 
$R$ of the form $H_{\alpha, m}$ separate $R$ from $A_\circ$, meaning that 
$(\alpha, x) > m$ holds for $x \in R$. The numbers $h_i (\Phi, m)$ were 
introduced and studied in \cite{Ath2}. Let $h_i (\Delta^m (\Phi))$ and $h_i 
(\Delta^m_+ (\Phi))$ be the $i$th entries of the $h$-vector of the simplicial 
complexes $\Delta^m (\Phi)$ and $\Delta^m_+ (\Phi)$, respectively. It follows
from case by case computations in \cite{Ath2, FR2, Tz} that $h_i (\Phi, m) 
= h_i (\Delta^m (\Phi))$ for all $i$ when $\Phi$ is of classical type in the 
Cartan-Killing classification.  

We define $h^+_i (\Phi, m)$ as the number of bounded regions $R$ of 
$\aA_\Phi^m$ in the fundamental chamber of $\aA_\Phi$ for which exactly 
$\ell - i$ walls of $R$ of the form $H_{\alpha, m}$ do not separate $R$ from 
the fundamental alcove $A_\circ$. Theorem \ref{thm0} implies that the sum 
of the numbers $h^+_i (\Phi, m)$, as well as that of $h_i (\Delta^m_+ 
(\Phi))$, for $0 \le i \le \ell$ is equal to $N^+ (\Phi, m)$. The 
significance of the numbers $h^+_i (\Phi, m)$ comes from the following 
conjecture, which can be viewed as the positive analogue of 
\cite[Conjecture 3.1]{FR2}. 

\begin{conjecture} 
For any irreducible crystallographic root system $\Phi$ and all $m \ge 1$ and 
$0 \le i \le \ell$ we have $h^+_i (\Phi, m) = h_i (\Delta^m_+ (\Phi))$.
\label{conj0}
\end{conjecture}

Our first main result (Corollary \ref{cor:conj}) establishes the previous 
conjecture when $m=1$ and when $\Phi$ has type $A$, $B$ or $C$ and $m$ is 
arbitrary. Our second main result provides 
combinatorial interpretations to the numbers $h^+_i (\Phi, m)$ similar to 
the ones given in \cite{Ath2} for $h_i (\Phi, m)$. To state this result we 
need to recall (or modify) some definitions
and notation from \cite{Ath2}. For $y \in T_m$ consider the stabilizer of 
$y$ with respect to the $W$-action on $T_m$. This is a subgroup of $W$ 
generated by reflections. The minimum number of reflections needed to 
generate this subgroup is its \emph{rank} and is denoted by $r(y)$. We 
may use the notation $r(x)$ for a $W$-orbit $x$ in $T_m$ since stabilizers 
of elements of $T_m$ in the same $W$-orbit are conjugate subgroups of $W$ 
and hence have the same rank. A subset $\jJ$ of $\Phi^+$ is an \emph{order 
ideal} if $\Phi^+ \, \sm \jJ$ is a filter. An increasing chain $\jJ_1 
\subseteq \jJ_2 \subseteq \cdots \subseteq \jJ_m$ of ideals in $\Phi^+$ is 
a \textit{geometric chain of ideals} of length $m$ if
\begin{equation}
(\jJ_i + \jJ_j) \, \cap \Phi^+ \subseteq \jJ_{i+j}
\label{bi2}
\end{equation}
holds for all indices $i, j$ with $i + j \le m$ and
\begin{equation}
(\iI_i + \iI_j) \, \cap \Phi^+ \subseteq \iI_{i+j}
\label{bi1}
\end{equation}
holds for all indices $i, j$, where $\iI_i = \Phi^+ \, \sm \jJ_i$ for 
$0 \le i \le m$ and $\iI_i = \iI_m$ for $i > m$. Such a chain is called 
\emph{positive} if $\jJ_m$ contains the set of simple roots or, 
equivalently, if $\iI_m$ is a positive filter. A positive root $\alpha$ 
is \emph{indecomposable} of \emph{rank} $m$ with respect to this 
increasing chain of ideals if $\alpha$ is a maximal element of $\jJ_m 
\sm \jJ_{m-1}$ and it is not possible to write $\alpha = \beta + \gamma$ 
with $\beta \in \jJ_i$ and $\gamma \in \jJ_j$ for indices $i, j \ge 1$ 
with $i + j = m$. The following theorem refines part of Theorem 
\ref{thm0}.  

\begin{theorem} 
Let $\Phi$ be an irreducible crystallographic root system of rank 
$\ell$ with Weyl group $W$, $m$ be a positive integer and $O_m (\Phi)$ 
be the set of orbits of the action of $W$ on $\check{Q} / \, (mh-1) \, 
\check{Q}$. For any $0 \le i \le \ell$ the following are equal:
\begin{enumerate}
\itemsep=0pt
\item[{\rm (i)}] the number $h^+_{\ell-i} (\Phi, m)$,
\item[{\rm (ii)}] the number of positive geometric chains of ideals in the 
root poset $\Phi^+$ of length $m$ having $i$ indecomposable elements 
of rank $m$,
\item[{\rm (iii)}] the number of orbits $x \in O_m (\Phi)$ with $r(x) = i$
and
\item[{\rm (iv)}] the number of points in $\check{Q} \cap (mh-1) \, 
\overline{A_\circ}$ which lie in $i$ walls of $(mh-1) \, 
\overline{A_\circ}$.
\end{enumerate}
In particular, the number of positive geometric chains of ideals in 
$\Phi^+$ of length $m$ is equal to $N^+ (\Phi, m)$.
\label{thm1}
\end{theorem}

The equivalence of (iii) and (iv) follows essentially from the 
results of \cite[Section 7.4]{Ha}. 
In the special case $m=1$ the arrangement $\aA_\Phi^m$ consists of 
the hyperplanes $H_\alpha$ and $H_{\alpha, 1}$ for all $\alpha \in 
\Phi$ and is known as the \textit{Catalan arrangement} associated to 
$\Phi$, denoted $\Cat_\Phi$. Moreover a geometric chain of ideals 
consists of a single ideal $\jJ$ in $\Phi$. This chain is positive if 
$\jJ$ contains the set of simple roots or, equivalently, if $\iI = 
\Phi^+ \, \sm \jJ$ is a positive filter and in that case the set 
of rank one indecomposable elements is the set of maximal elements 
of $\jJ$. We write $h^+_i (\Phi)$ instead of $h^+_i (\Phi, m)$ when 
$m=1$. Part of the next corollary is implicit in the work of 
E.~Sommers \cite[Section 6]{So}.

\begin{corollary}
Let $\Phi$ be an irreducible crystallographic root system of rank 
$\ell$ with Weyl group $W$ and $O (\Phi)$ be the set of orbits of the 
action of $W$ on $\check{Q} / (h-1) \, \check{Q}$. For any $0 \le i \le 
\ell$ the following are equal:
\begin{enumerate}
\itemsep=0pt
\item[{\rm (i)}] the number of ideals in the root poset $\Phi^+$ which 
contain the set of simple roots and have $i$ maximal elements,
\item[{\rm (ii)}] the number $h^+_{\ell-i} (\Phi)$ of bounded regions $R$ of 
$\Cat_\Phi$ in the fundamental chamber of $\aA_\Phi$ such that $i$ walls 
of $R$ of the form $H_{\alpha, 1}$ do not separate $R$ from $A_\circ$,
\item[{\rm (iii)}] the number of orbits $x \in O (\Phi)$ with $r(x) = i$,
\item[{\rm (iv)}] the number of points in $\check{Q} \cap (h-1) \, 
\overline{A_\circ}$ which lie in $i$ walls of $(h-1) \, 
\overline{A_\circ}$ and
\item[{\rm (v)}] the entry $h_{\ell-i} (\Delta_+ (\Phi))$ of the $h$-vector of
the positive part of $\Delta (\Phi)$.
\end{enumerate}
\label{cor1}
\end{corollary}

Our last theorem provides a different interpretation to the numbers $h^+_i 
(\Phi)$ in terms of order filters in $\Phi^+$.
\begin{theorem} 
For any irreducible crystallographic root system $\Phi$ and any nonnegative
integer $i$ the number $h^+_i (\Phi)$ is equal to the number of positive 
filters in $\Phi^+$ having $i$ minimal elements.
\label{thm2}
\end{theorem}

\vspace{0.1 in}
Theorem \ref{thm1} is proved in Sections \ref{som} and \ref{proof} by 
means of two bijections. The first is the restriction of a bijection 
of \cite[Section 3]{Ath2} and maps the set of positive geometric 
chains of ideals in $\Phi^+$ of length $m$ to the set of bounded regions 
of $\aA_\Phi^m$ in the fundamental chamber (Theorem \ref{cor:bij1}) while
the second maps this set of regions to the set of $W$-orbits of $T_m$.
In the case $m=1$ the composite of these two bijections gives essentially
a bijection of Sommers \cite{So} from the set of positive filters in 
$\Phi^+$ to $\check{Q} \cap (h-1) \, \overline{A_\circ}$. The proof of
Theorem \ref{thm1} in these two sections parallels the one of Theorem 1.2 
in \cite{Ath2} and for this reason most of the details are omitted. The
main difference is that the unique alcove in a fixed bounded region of 
$\aA_\Phi^m$ which is furthest away from $A_\circ$ plays the role played 
in \cite{Ath2} by the unique alcove in a region of $\aA_\Phi^m$ closest to 
$A_\circ$. The existence of these maximal alcoves was first established 
and exploited in the special case $m=1$ by Sommers \cite{So}. In Section 
\ref{f} we prove Conjecture \ref{conj0} when $m=1$ and when $\Phi$ has type 
$A$, $B$ or $C$ and $m$ is arbitrary (Corollary \ref{cor:conj}) using the 
fact that $h_i (\Phi, m) = h_i (\Delta^m (\Phi))$ holds for all $i$ in 
these cases. A key ingredient in the proof is a new combinatorial 
interpretation (see part (iii) of Theorems \ref{thm:eleni1} and 
\ref{thm:eleni2}) to the $f$-numbers defined from the $h_i (\Phi, m)$ 
and $h^+_i (\Phi, m)$ via the usual identity relating $f$-vectors and 
$h$-vectors of simplicial complexes. In Section \ref{class} we compute 
the numbers which appear in Theorem \ref{thm1} for root systems of 
classical type and those in Corollary \ref{cor1} for root systems of 
exceptional type. We also prove Theorem \ref{thm2} by exploiting the 
symmetry of the distribution of the set of all filters in $\Phi^+$ by the 
number of minimal elements, observed by D. Panyushev \cite{Pa}. Some useful 
background material is summarized in Section \ref{pre}. We conclude with 
some remarks in Section 
\ref{remarks}.

Apart from \cite{Ath2}, our motivation for this work comes to a great 
extent from the papers by Fomin and Reading \cite{FR2}, Fomin and 
Zelevinsky \cite{FZ} and Sommers \cite{So}.

\section{Preliminaries}
\label{pre}

In this section we introduce notation and terminology and recall a few 
useful facts related to root systems, affine Weyl groups, generalized
cluster complexes and the combinatorics of $\aA^m_\Phi$. We refer to 
\cite{Hu} and \cite{Ath2, FR1, FR2} for further background and 
references and warn the reader that, throughout the paper, some of our 
notation and terminology differs from that employed in \cite{Ath2} (this
is done in part to ease the co-existence of order filters and order 
ideals in this paper, typically denoted by the letters $\iI$ and $\jJ$, 
respectively, and in part to match some of the notation of \cite{FR2}).  

\medskip
\noindent
{\bf Root systems and Weyl groups.}
Let $V$ be an $\ell$-dimensional Euclidean space with inner product 
$( \ , \ )$. Given a hyperplane arrangement $\aA$ in $V$, meaning a 
discrete set of affine subspaces of $V$ of codimension one, the 
\emph{regions} of $\aA$ are the connected components of the space 
obtained from $V$ by removing the hyperplanes in $\aA$. Let $\Phi$ be 
a crystallographic root system spanning $V$. For any real $k$ and 
$\alpha \in \Phi$ we denote by $H_{\alpha, k}$ the hyperplane in $V$ 
defined by the equation $(\alpha,x) = k$ and set $H_{\alpha} = H_{\alpha, 
0}$. We fix a positive system $\Phi^+ \subseteq \Phi$ and the 
corresponding (ordered) set of simple roots $\Pi = 
\{\sigma_1,\ldots,\sigma_\ell\}$. For $1 \le i \le \ell$ we denote by 
$s_i$ the orthogonal reflection in the hyperpane $H_{\sigma_i}$, called a 
\emph{simple reflection}. We will often write $\Phi_I$ 
instead of $\Phi$, where $I$ is an index set in bijection with $\Pi$, 
and denote by $\Phi_J$ the parabolic root system corresponding to $J 
\subseteq I$. If $\Phi$ is irreducible we denote by $\tilde{\alpha}$ 
the highest positive root, by $e_1, e_2,\dots,e_{\ell}$ the exponents 
and by $h$ the Coxeter number of $\Phi$ and set $p=mh-1$, where $m$ is 
a fixed positive integer. The following well known lemmas will be used, 
as in \cite{Ath2}.
\begin{lemma} {\rm (\cite[Lemma 2.1]{Ath2})}
{\rm (i)} If $\alpha_1, \alpha_2,\dots,\alpha_r \in \Phi^+$ with
$r \ge 2$ and $\alpha = \alpha_1 + \alpha_2 + \cdots + \alpha_r \in 
\Phi^+$ then there exists $i$ with $1 \le i \le r$ such that $\alpha 
- \alpha_i \in \Phi^+$.

\noindent
{\rm (ii)} {\rm (cf. \cite{Pa, So})} 
If $\alpha_1, \alpha_2,\dots,\alpha_r \in \Phi$ and 
$\alpha_1 + \alpha_2 + \cdots + \alpha_r = \alpha \in \Phi$ then 
$\alpha_1 = \alpha$ or there exists $i$ with $2 \le i \le r$ such 
that $\alpha_1 + \alpha_i \in \Phi \cup \{0\}$. \qed
\label{lem:ro}
\end{lemma}
\begin{lemma} {\rm (\cite[Ch. 6, 1.11, Proposition 31]{Bou} \cite[p. 
84]{Hu})} If $\Phi$ is irreducible and $\tilde{\alpha} = \sum_{i=1}^{\ell} 
c_i \, \sigma_i$ then $\sum_{i=1}^{\ell} c_i = h-1$.
\label{lem:h}
\end{lemma}
We denote by $\aA_{\Phi}$ the 
\emph{Coxeter arrangement} associated to $\Phi$, i.e. the collection 
of linear hyperplanes $H_\alpha$ in $V$ with $\alpha \in \Phi$, and 
by $W$ the corresponding \emph{Weyl group}, generated by the reflections 
in these hyperplanes. Thus $W$ is finite and minimally generated by the 
set of simple reflections, it leaves $\Phi$ invariant and acts simply 
transitively on the set of regions of $\aA_\Phi$, called \emph{chambers}. 
The \emph{fundamental chamber} is the region defined by the inequalities 
$0 < (\alpha, x)$ for $\alpha \in \Phi^+$. A subset of $V$ is called 
\emph{dominant} if it is contained in the fundamental chamber. The 
\emph{coroot lattice} $\check{Q}$ of $\Phi$ is the $\ZZ$-span of the 
set of coroots
\[ \Phi^\vee = \left\{ \frac{2 \alpha}{(\alpha,\alpha)}: \ \alpha \in 
\Phi \right\}. \] 

>From now on we assume for simplicity that $\Phi$ is irreducible. The 
group $W$ acts on the lattice $\check{Q}$ and on its sublattice $p \, 
\check{Q}$, hence it also acts on the quotient $T_m (\Phi) = \check{Q} 
/ p \, \check{Q}$. We denote by $O_m (\Phi)$ the set of orbits of the 
$W$-action on $T_m (\Phi)$ and use the notation $T (\Phi)$ and $O 
(\Phi)$ when $m=1$. We denote 
by $\widetilde{\aA}_\Phi$ the \emph{affine Coxeter arrangement}, 
which is the infinite hyperplane arrangement in $V$ consisting of the 
hyperplanes $H_{\alpha, k}$ for $\alpha \in \Phi$ and $k \in \ZZ$, and 
by $W_a$ the \emph{affine Weyl group}, generated by the reflections in 
the hyperplanes of $\widetilde{\aA}_\Phi$. The group $W_a$ is the 
semidirect product of $W$ and the translation group in $V$ corresponding 
to the coroot lattice $\check{Q}$ and is minimally generated by the set 
$\{s_0, s_1,\dots,s_{\ell}\}$ of \emph{simple affine reflections}, where 
$s_0$ is the reflection in the hyperplane $H_{\widetilde{\alpha}, 1}$. 
For $w \in W_a$ and $0 \le i \le \ell$, the reflection $s_i$ is a 
\emph{right ascent} of $w$ if $\ell (ws_i) > \ell(w)$, where $\ell (w)$ 
is the length of the shortest expression of $w$ as a product of simple 
affine reflections. 
The group $W_a$ acts simply transitively on the set of regions of 
$\widetilde{\aA}_{\Phi}$, called \emph{alcoves}. The \emph{fundamental 
alcove} of $\widetilde{\aA}_{\Phi}$ can be defined as
\[ A_\circ = \{ x \in V: \, 0 < (\sigma_i, x ) \ 
{\rm for} \ 1 \le i \le \ell \ {\rm and} \ (\tilde{\alpha}, x) < 
1\}. \]
Note that every alcove can be written as $w A_{\circ}$ for a unique $w 
\in W_a$. Moreover, given $\alpha \in \Phi^+$, there exists a unique
integer $r$, denoted $r(w, \alpha)$, such that $r-1 < (\alpha, x) < r$
holds for all $x \in w A_\circ$. The next lemma is a reformulation of
the main result of \cite{Sh0}.
\begin{lemma} \mbox~{\rm (\cite[Theorem 5.2]{Sh0}).}
Let $r_\alpha$ be an integer for each $\alpha \in \Phi^+$. There exists
$w \in W_a$ such that $r(w, \alpha) = r_\alpha$ for each $\alpha \in 
\Phi^+$ if and only if
\begin{equation}
r_\alpha + r_\beta - 1 \le r_{\alpha + \beta} \le r_\alpha + r_\beta
\label{eq:alc}
\end{equation}
for all $\alpha, \beta \in \Phi^+$ with $\alpha + \beta \in \Phi^+$.
\label{lem:alc}
\end{lemma}
We say that two open regions in $V$ are \emph{separated} by 
a hyperplane $H \in \widetilde{\aA}_{\Phi}$ if they lie in different 
half-spaces relative to $H$. If $R$ is a region of a subarrangement of 
$\widetilde{\aA}_{\Phi}$ or the closure of such a region (in particular, 
if $R$ is a chamber or an alcove), we refer to the hyperplanes of 
$\widetilde{\aA}_\Phi$ which support facets of the closure of $R$ as 
the \emph{walls} of $R$.

\medskip
\noindent
{\bf Generalized cluster complexes.} Let $\Phi$ be crystallographic 
(possibly reducible) of rank $\ell$. The generalized cluster complex 
$\Delta^m (\Phi)$ is an abstract simplicial complex on the vertex set 
$\Phi^m_{\ge -1}$ consisting of the negative simple roots and $m$ copies 
of each positive root; we refer to \cite[Section 1.2]{FR2} for the 
definition. It is a pure complex of dimension $\ell-1$ \cite[Proposition 
1.7]{FR2}. If $\Phi$ is a direct product $\Phi = \Phi_1 \times \Phi_2$ 
then $\Delta^m (\Phi)$ is the simplicial join of $\Delta^m (\Phi_1)$ and 
$\Delta^m (\Phi_2)$. We denote by $\Delta^m_+ (\Phi)$ the induced 
subcomplex of $\Delta^m (\Phi)$ on the set of vertices obtained from 
$\Phi^m_{\ge -1}$ by removing the negative simple roots and call this 
simplicial complex the \emph{positive part} of $\Delta^m (\Phi)$. For $0 
\le i \le \ell$ we denote by $f_{i-1} (\Delta^m (\Phi))$ and 
$f_{i-1} (\Delta^m_+ (\Phi))$ the number of $(i-1)$-dimensional faces of 
the complex $\Delta^m (\Phi)$ and $\Delta^m_+ (\Phi)$, respectively. These 
numbers are related to the $h_i (\Delta^m (\Phi))$ and $h_i (\Delta^m_+ 
(\Phi))$ by the equations
\begin{equation}
\sum_{i=0}^\ell \, f_{i-1} (\Delta^m (\Phi)) (x-1)^{\ell-i} = \sum_{i=0}^\ell 
\, h_i (\Delta^m (\Phi)) \, x^{\ell-i} 
\label{fFR}
\end{equation}
and
\begin{equation}
\sum_{i=0}^\ell \, f_{i-1} (\Delta^m_+ (\Phi)) (x-1)^{\ell-i} = 
\sum_{i=0}^\ell \, h_i (\Delta^m_+ (\Phi)) \, x^{\ell-i} 
\label{f+FR}
\end{equation}
respectively. 

Following \cite{FR2, Tz} we give explicit combinatorial descriptions of the
complexes $\Delta^m (\Phi)$ and $\Delta^m_+ (\Phi)$ when $\Phi$ has type $A$,
$B$ or $C$. For $\Phi = A_{n-1}$ let $\ppP$ be a convex polygon with $mn+2$
vertices. A diagonal of $\ppP$ is called \emph{$m$-allowable} if it divides
$\ppP$ into two polygons each with number of vertices congruent to 2 mod 
$m$. Vertices of $\Delta^m (\Phi)$ are the $m$-allowable diagonals of $\ppP$
and faces are the sets of pairwise noncrossing diagonals of this kind. For
$\Phi = B_n$ or $C_n$ let $\qqQ$ be a centrally symmetric convex polygon with 
$2mn+2$ vertices. A vertex of $\Delta^m (\Phi)$ is either a diameter of 
$\qqQ$, i.e. a diagonal connecting antipodal vertices, or a pair of 
$m$-allowable diagonals related by a half-turn about the center of $\qqQ$. 
A set of vertices of $\Delta^m (\Phi)$ forms a face if the diagonals of 
$\qqQ$ defining these vertices are pairwise noncrossing. In all cases the 
explicit bijection of $\Phi^m_{\ge -1}$ with the set of allowable diagonals 
of $\ppP$ or $\qqQ$ just described is analogous to the one given in 
\cite[Section 3.5]{FZ} for the usual cluster complex $\Delta (\Phi)$, so 
that the negative simple roots form an $m$-snake of allowable diagonals in 
$\ppP$ or $\qqQ$ and $\Delta^m_+ (\Phi)$ is the subcomplex of $\Delta^m 
(\Phi)$ obtained by removing the vertices in the $m$-snake.

The \emph{negative part} of a face $c$ of $\Delta^m (\Phi)$
is the set of indices $J \subseteq I$, where $\Phi = \Phi_I$, which correspond
to the negative simple roots contained in $c$. The next lemma appears as 
\cite[Proposition 3.6]{FZ} in the case $m=1$ and follows from the explicit 
description of the relevant complexes in the remaining cases.
\begin{lemma}
Assume that either $m=1$ or $\Phi_I$ has type $A, B$ or $C$. For any $J 
\subseteq I$ the map $c \mapsto c \sm \, \{-\alpha_i: i \in J\}$ 
is a bijection from the set of faces of $\Delta^m (\Phi_I)$ with negative 
part $J$ to the set of faces of $\Delta^m_+ (\Phi_{I \sm J})$. In particular
\begin{equation}
f_{k-1} (\Delta^m (\Phi_I)) = \sum_{J \subseteq I} \ f_{k-|J|-1} 
(\Delta^m_+ (\Phi_{I \sm J})),
\label{ff+FR}
\end{equation}
where $f_{i-1} (\Delta) = 0$ if $i < 0$ for any complex $\Delta$. \qed
\label{lem:FR}
\end{lemma}

For $m=1$ essentially the same equation as (\ref{ff+FR}) has appeared in the 
context of quiver representations in \cite[Section 6]{MRZ}. 

\medskip
\noindent
{\bf Regions of $\aA^m_\Phi$ and chains of filters.} Let $\Phi$ be 
irreducible and crystallographic of rank $\ell$. For $0 \le i \le \ell$ let 
$h_i (\Phi, m)$ be the number of dominant regions $R$ of $\aA_\Phi^m$ for 
which exactly $\ell - i$ walls of $R$ of the form $H_{\alpha, m}$ separate 
$R$ from $A_\circ$, as in Section \ref{intro}. We recall another combinatorial 
interpretation of $h_i (\Phi, m)$ from \cite{Ath2} using slightly different 
terminology. We call a decreasing chain 
$$\Phi^+ = \iI_0 \supseteq \iI_1 \supseteq \iI_2 \supseteq \cdots \supseteq 
\iI_m$$ 
of filters in $\Phi^+$ 
a \textit{geometric chain of filters} of length $m$ if (\ref{bi2}) and 
(\ref{bi1}) hold under the same conventions as in Section \ref{intro} (the 
term \textit{co-filtered chain of dual order ideals} was used in 
\cite{Ath2} instead). A positive root $\alpha$ is \emph{indecomposable} of 
\emph{rank} $m$ with respect to this chain if $\alpha \in \iI_m$ and it is 
not possible to write $\alpha = \beta + \gamma$ with $\beta \in \iI_i$ and 
$\gamma \in \iI_j$ for indices $i, j \ge 0$ with $i + j = m$. Let $R_\iI$ 
be the set of points $x \in V$ which satisfy
\begin{equation}
\begin{tabular}{ll} $(\alpha, x) > r$, & {\rm if} \ $\alpha \in \iI_r$
\\ $0 < (\alpha, x) < r$, & {\rm if} \ $\alpha \in \jJ_r$
\end{tabular} 
\label{map}
\end{equation}
for $0 \le r \le m$, where $\jJ_r = \Phi^+ \, \sm \iI_r$. The following 
statement combines parts of Theorems 3.6 and 3.11 in \cite{Ath2}.  

\begin{theorem}
The map $\iI \mapsto R_\iI$ is a bijection from the set of geometric chains 
of filters of length $m$ in $\Phi^+$ to the set of dominant regions of 
$\aA^m_\Phi$. Moreover a positive root $\alpha$ is indecomposable of rank $m$ 
with respect to $\iI$ if and only if $H_{\alpha, m}$ is a wall of $R_\iI$
which separates $R_\iI$ from $A_\circ$.
 
In particular $h_i (\Phi, m)$ is equal to the number of geometric chains 
of filters in $\Phi^+$ of length $m$ having $\ell-i$ indecomposable elements 
of rank $m$. \hfill
\qed
\label{thm:ath2}
\end{theorem}

By modifying the definition given earlier or using the interpretation in 
the last statement of the previous theorem we can define the numbers $h_i 
(\Phi, m)$ when $\Phi$ is reducible as well. Clearly 
\[ h_k (\Phi_1 \times \Phi_2, m) = \sum_{i + j = k} \ h_i (\Phi_1, 
m) \, h_j (\Phi_2, m) \]
for any crystallographic root systems $\Phi_1, \Phi_2$.

\section{Chains of ideals, bounded begions and maximal alcoves}
\label{som}

In this section we generalize some of the results of Sommers \cite{So} on 
bounded dominant regions of $\Cat_\Phi$ and positive filters in $\Phi^+$ to 
bounded dominant regions of $\aA_\Phi^m$ and positive geometric chains of 
ideals and establish the equality of the numbers appearing in (i) and (ii) 
in the statement of Theorem \ref{thm1}. The results of this and the following 
section are analogues of the results of Sections 3 and 4 of \cite{Ath2} 
on the set of all dominant regions of $\aA_\Phi^m$. Their proofs are 
obtained by minor adjustments from those of \cite{Ath2}, suggested by the 
modifications of the relevant definitions, and thus are only sketched or 
omitted.

Let $\Phi$ be irreducible and crystallographic of rank $\ell$ and let 
$\jJ$ be a positive geometric chain of ideals 
\[ \emptyset = \jJ_0 \subseteq \jJ_1 \subseteq \jJ_2 \subseteq \cdots 
\subseteq \jJ_m \]
in $\Phi^+$ of length $m$, so that (\ref{bi2}) and (\ref{bi1}) hold, 
where $\iI_i = \Phi^+ \sm \, \jJ_i$, and $\Pi \subseteq \jJ_m$. We 
define
\[ r_\alpha(\jJ) = \min \{r_1 + r_2 + \cdots + r_k: \alpha = \alpha_1 
+ \alpha_2 + \cdots + \alpha_k \ {\rm with} \ \alpha_i \in \jJ_{r_i} \ 
{\rm for \ all} \ i\} \]
for any $\alpha \in \Phi^+$. Observe that $r_\alpha(\jJ)$ is well defined
since $\Pi \subseteq \jJ_m$ and that $r_\alpha (\jJ) \le r$ for 
$\alpha \in \jJ_r$, with $r_\alpha (\jJ) = 1$ if and only if $\alpha 
\in \jJ_1$. 
\begin{lemma}
If $\alpha = \alpha_1 + \alpha_2 + \cdots + \alpha_k \in \Phi^+$ and 
$\alpha_i \in \Phi^+$ for all $i$ then 
\[ r_\alpha(\jJ) \, \le \, \sum_{i=1}^k \, r_{\alpha_i} (\jJ). \]
\label{lem0}
\end{lemma}

\noindent
\emph{Proof.}
This is clear from the definition.
\qed

\begin{lemma}
Let $\alpha \in \Phi^+$ and $r_\alpha(\jJ) = r$. 
\begin{enumerate}
\itemsep=0pt
\item[{\rm (i)}] If $r \le m$ then $\alpha \in \jJ_r$.
\item[{\rm (ii)}] If $r > m$ then there exist $\beta, \gamma \in \Phi^+$ with 
$\alpha = \beta + \gamma$ and $r = r_\beta(\jJ) + r_\gamma(\jJ)$. Moreover 
we may choose $\beta$ so that $r_\beta (\jJ) \le m$.
\end{enumerate}
\label{lem1}
\end{lemma}

\noindent
\emph{Proof.}
Analogous to the proof of \cite[Lemma 3.2]{Ath2}.
\qed

\begin{lemma}
If $\alpha, \beta, \alpha + \beta \in \Phi^+$ and $a, b$ are 
integers such that $r_{\alpha + \beta}(\jJ) \le a+b$ then 
$r_\alpha(\jJ) \le a$ or $r_\beta(\jJ) \le b$.
\label{lem:cp}
\end{lemma}

\noindent
\emph{Proof.}
By induction on $r_{\alpha + \beta}(\jJ)$, as in the proof of \cite[Lemma 
3.3]{Ath2}.
\qed

\begin{corollary}
We have 
\[ r_\alpha(\jJ) + r_\beta(\jJ) - 1 \, \le \, r_{\alpha + \beta}(\jJ) \, 
\le \, r_\alpha(\jJ) + r_\beta(\jJ) \]
whenever $\alpha, \beta, \alpha + \beta \in \Phi^+$. 
\label{cor:k}
\end{corollary}

\noindent
\emph{Proof.}
The second inequality is a special case of Lemma \ref{lem0} and the first
follows from Lemma \ref{lem:cp} letting $a = r_\alpha (\jJ) - 1$ and $b = 
r_{\alpha + \beta}(\jJ) - a$.
\qed

\medskip
We denote by $R_\jJ$ the set of points $x \in V$ which satisfy the 
inequalities in (\ref{map}) (thus we allow a slight abuse of notation 
since the same set was denoted by $R_\iI$ in Section \ref{pre}, where 
$\iI$ is the chain of complementary filters $\iI_i$). Since $\Pi \subseteq 
\jJ_m$ we have $0 < (\sigma_i, x) < m$ for all $1 \le i \le \ell$ and 
$x \in R_\jJ$ and therefore $R_\jJ$ is bounded.
\begin{proposition}
There exists a unique $w \in W_a$ such that $r(w, \alpha) = r_\alpha(\jJ)$ 
for $\alpha \in \Phi^+$. Moreover, $w A_\circ \subseteq R_\jJ$. In particular, 
$R_\jJ$ is nonempty.
\label{prop:min}
\end{proposition}

\noindent
\emph{Proof.}
The existence in the first statement follows from Lemma \ref{lem:alc} and 
Corollary \ref{cor:k} while uniqueness is obvious.
For the second statement let $\alpha \in \Phi^+$ and $1 \le r \le m$. Part (i)
of Lemma \ref{lem1} implies that $r_\alpha (\jJ) \le r$ if and only if $\alpha 
\in J_r$. Hence from the inequalities
\[ r_\alpha (\jJ) - 1 \, < \, (\alpha, x) \, < \, r_\alpha (\jJ), \]
which hold for $x \in w A_\circ$, we conclude that $w A_\circ \subseteq 
R_\jJ$. 
\qed

\medskip
Let $\psi$ be the map which assigns the set $R_\jJ$ to a positive geometric
chain of ideals $\jJ$ in $\Phi^+$ of length $m$. Conversely, given a bounded
dominant region $R$ of $\aA_\Phi^m$ let 
$\phi(R)$ be the sequence $\emptyset = \jJ_0 \subseteq \jJ_1 \subseteq \jJ_2 
\subseteq \cdots \subseteq \jJ_m$ where $\jJ_r$ is the set of $\alpha \in 
\Phi^+$ for which $(\alpha, x) < r$ holds in $R$. Clearly each $\jJ_r$ is an
ideal in $\Phi^+$.
\begin{theorem}
The map $\psi$ is a bijection from the set of positive geometric chains of 
ideals in $\Phi^+$ of length $m$ to the set of bounded dominant regions of 
$\aA_\Phi^m$, and the map $\phi$ is its inverse.
\label{cor:bij1}
\end{theorem}

\noindent
\emph{Proof.}
That $\psi$ is well defined follows from Proposition \ref{prop:min}, 
which guarantees that $R_\jJ$ is nonempty (and bounded). To check that 
$\phi$ is well defined observe that if $R$ is a bounded dominant region 
of $\aA_\Phi^m$ and if $(\alpha, x) < i$ and $(\beta, 
x) < j$ hold for $x \in R$ then $(\alpha + \beta, x) < i+j$ must hold for 
$x \in R$, so that $\phi(\jJ)$ satisfies (\ref{bi2}). Similarly, $\phi(\jJ)$ 
satisfies (\ref{bi1}). That $\Pi \subseteq \jJ_m$ follows from \cite[Lemma 
4.1]{Ath1}. It is clear that $\psi$ and $\phi$ are inverses of each other. 
\qed

\medskip
Let $R = R_\jJ$ be a bounded dominant region of $\aA_\Phi^m$, where $\jJ 
= \phi(R)$. Let $w_R$ denote the element of the affine Weyl group $W_a$ 
which is assigned to $\jJ$ in Proposition \ref{prop:min}. The following 
proposition implies that $w_R A_\circ$ is the alcove in $R$ which is the 
furthest away from $A_\circ$. In the special case $m=1$ the existence 
of such an alcove was established by Sommers \cite[Proposition 5.4]{So}.


%
\begin{proposition}
Let $R$ be a bounded dominant region of $\aA_\Phi^m$. The element $w_R$ 
is the unique $w \in W_a$ such that 
$w A_\circ \subseteq R$ and whenever $\alpha \in \Phi^+$, $r \in \ZZ$ 
and $(\alpha, x) > r$ holds for some $x \in R$ we have $(\alpha, x) 
> r$ for all $x \in w A_\circ$. 
\label{prop:wR}
\end{proposition}

\noindent
\emph{Proof.}
Analogous to the proof of \cite[Proposition 3.7]{Ath2}.
\qed

\medskip
We now introduce the notion of an indecomposable element with respect 
to the increasing chain of ideals $\jJ$.
\begin{definition} 
Given $1 \le r \le m$, a root $\alpha \in \Phi^+$ is indecomposable of 
rank $r$ with respect to $\jJ$ if $\alpha \in \jJ_r$ and
\begin{enumerate}
\itemsep=0pt
\item[{\rm (i)}] $r_\alpha(\jJ) = r$,
\item[{\rm (ii)}] it is not possible to write $\alpha = \beta + \gamma$ with 
$\beta \in \jJ_i$ and $\gamma \in \jJ_j$ for indices $i, j \ge 1$ with 
$i + j = r$ and
\item[{\rm (iii)}] if $r_{\alpha + \beta} (\jJ) = t \le m$ for some $\beta 
\in \Phi^+$ then $\beta \in \jJ_{t-r}$.    
\end{enumerate}
\label{def}
\end{definition}

Observe that, by part (i) of Lemma \ref{lem1}, the assumption $\alpha 
\in \jJ_r$ in this definition is actually implied by condition (i). For
$r=m$ the definition is equivalent to the one proposed in Section 
\ref{intro}, as the following lemma shows.
\begin{lemma}
A positive root $\alpha$ is indecomposable of rank $m$ with respect to 
$\jJ$ if and only if $\alpha$ is a maximal element of $\jJ_m \sm \jJ_{m-1}$ 
and it is not possible to write $\alpha = \beta + \gamma$ with $\beta \in 
\jJ_i$ and $\gamma \in \jJ_j$ for indices $i, j \ge 1$ with $i + j = m$.
\label{lem:m}
\end{lemma}

\noindent
\emph{Proof.}
Suppose that $\alpha \in \jJ_m$ is indecomposable of rank $m$. Since 
$r_\alpha(\jJ) = m$ we must have $\alpha \notin \jJ_{m-1}$. Hence to show 
that $\alpha$ satisfies the condition in the statement of the lemma
it suffices to show that $\alpha$ is maximal in $\jJ_m$. If not then by 
Lemma \ref{lem:ro} (i) there exists $\beta \in \Phi^+$ such that $\alpha 
+ \beta \in \jJ_m$. Then clearly $r_{\alpha + \beta} (\jJ) \le m$ and  
$r_{\alpha + \beta} (\jJ) \ge m$ by Corollary \ref{cor:k}. Hence 
$r_{\alpha + \beta} (\jJ) = m$  and condition (iii) of Definition 
\ref{def} leads to a contradiction.
 
For the converse, suppose that $\alpha \in \jJ_m$ satisfies the condition 
in the statement of the lemma. In view of part (i) of Lemma \ref{lem1},
condition (iii) in Definition \ref{def} 
is satisfied since $\alpha$ is assumed to be maximal in $\jJ_m$. Hence to 
show that $\alpha$ is indecomposable of rank $m$ it suffices to show that 
$r_\alpha(\jJ) = m$. This is implied by the assumption that $\alpha \notin
\jJ_{m-1}$ and part (i) of Lemma \ref{lem1}. 
\qed

\begin{lemma}
Suppose that $\alpha$ is indecomposable with respect to $\jJ$.
\begin{enumerate}
\itemsep=0pt
\item[{\rm (i)}] We have $r_\alpha(\jJ) = r_\beta(\jJ) + r_\gamma(\jJ) - 1$
whenever $\alpha = \beta + \gamma$ with $\beta, \gamma \in \Phi^+$.
\item[{\rm (ii)}] We have $r_\alpha(\jJ) + r_\beta(\jJ) = r_{\alpha + 
\beta}(\jJ)$ whenever $\beta, \alpha + \beta \in \Phi^+$.    
\end{enumerate}
\label{lem:ind}
\end{lemma}

\noindent
\emph{Proof.}
Analogous to the proof of \cite[Lemma 3.10]{Ath2}. For part (ii), letting
$r_\alpha(\jJ) = r$ and $r_{\alpha + \beta}(\jJ) = t$, we prove instead 
that $r_\beta(\jJ) \le t-r$. This implies the result by Corollary 
\ref{cor:k}. 
\qed

\medskip
The following theorem explains the connection between indecomposable 
elements of $\jJ$ and walls of $R_\jJ$. 
\begin{theorem}
If $\jJ$ is a positive geometric chain of ideals in $\Phi^+$ of length 
$m$ with corresponding region $R = R_\jJ$ and $1 \le r \le m$ then the 
following sets are equal:
\begin{enumerate}
\itemsep=0pt
\item[{\rm (i)}] the set of indecomposable roots $\alpha \in \Phi^+$ with 
respect to $\jJ$ of rank $r$, 
\item[{\rm (ii)}] the set of $\alpha \in \Phi^+$ such that $H_{\alpha, r}$ 
is a wall of $R$ which does not separate $R$ from $A_\circ$ and
\item[{\rm (iii)}] the set of $\alpha \in \Phi^+$ such that $H_{\alpha, r}$ 
is a wall of $w_R A_\circ$ which does not separate $w_R A_\circ$ from 
$A_\circ$.    
\end{enumerate}
\label{mylemma}
\end{theorem}

\noindent
\emph{Proof.}
We prove that $F_r (R) \subseteq F_r (\jJ) \subseteq F_r (w_R) \subseteq 
F_r (R)$ for the three sets defined in the statement of the theorem as 
in the proof of \cite[Theorem 3.11]{Ath2}, replacing the inequalities 
$(\alpha, x ) > k$ which appear there by $(\alpha, x ) < r$ and recalling 
from the proof of Proposition \ref{prop:wR} that $(\alpha, x) < \, 
r_\alpha (\jJ)$ holds for all $\alpha \in \Phi^+$ and $x \in R_\jJ$.
\qed

\medskip
We denote by $W_m (\Phi)$ the subset of $W_a$ consisting of the 
elements $w_R$ for the bounded dominant regions $R$ of $\aA_\Phi^m$; 
see Figure \ref{figure1} for the case $\Phi = A_2$ and $m=2$. We
abbreviate this set as $W (\Phi)$ in the case $m=1$. The elements of 
$W (\Phi)$ are called \emph{maximal} in \cite{So}.
\begin{corollary}
For any nonnegative integers $i_1, i_2,\dots,i_m$ the following
are equal:
\begin{enumerate}
\itemsep=0pt
\item[{\rm (i)}] the number of positive geometric chains of ideals in $\Phi^+$ of 
length $m$ having $i_r$ indecomposable elements of rank $r$ for each 
$1 \le r \le m$,
\item[{\rm (ii)}] the number of bounded dominant regions $R$ of $\aA_\Phi^m$ 
such that $i_r$ walls of $R$ of the form $H_{\alpha, r}$ do not separate 
$R$ from $A_\circ$ for each $1 \le r \le m$ and
\item[{\rm (iii)}] the number of $w \in W_m (\Phi)$ such that $i_r$ walls of 
$w A_\circ$ of the form $H_{\alpha, r}$ do not separate $w A_\circ$ from 
$A_\circ$ for each $1 \le r \le m$.
\end{enumerate}
\label{cor:ik}
\end{corollary}

\noindent
\emph{Proof.}
Combine Theorems \ref{cor:bij1} and \ref{mylemma}.
\qed

\medskip
The following corollary is immediate.
\begin{corollary}
For any nonnegative integer $i$ the numbers which appear in {\rm (i)} 
and {\rm (ii)} in the statement of Theorem \ref{thm1} are both equal to
the number of $w \in W_m (\Phi)$ such that $i$ walls of $w A_\circ$ of 
the form $H_{\alpha, m}$ do not separate $w A_\circ$ from $A_\circ$. 
\qed
\label{cor:part1}
\end{corollary}

As was the case with $h_i (\Phi, m)$, the interpretation in part (ii)
of Theorem \ref{thm1} mentioned in the previous corollary or the original
definition can be used to define $h^+_i (\Phi, m)$ when $\Phi$ is 
reducible. Equivalently we define
\[ h^+_k (\Phi_1 \times \Phi_2, m) = \sum_{i + j = k} \ h^+_i (\Phi_1, 
m) \, h^+_j (\Phi_2, m) \]
for any crystallographic root systems $\Phi_1, \Phi_2$.
We now consider the special case $m=1$. A positive geometric chain of 
ideals $\jJ$ of length $m$ in this case is simply a single ideal $\jJ$ in 
$\Phi^+$ such that $\Pi \subseteq \jJ$, meaning that $\iI = \Phi^+ \, \sm 
\jJ$ is a positive filter. By Lemma \ref{lem:m} the rank one indecomposable 
elements of $\jJ$ are exactly the maximal elements of $\jJ$. 
\begin{corollary}
For any nonnegative integer $i$ the following are equal to $h^+_{\ell-i} 
(\Phi)$:
\begin{enumerate}
\itemsep=0pt
\item[{\rm (i)}] the number of ideals in the root poset $\Phi^+$ which contain
all simple roots and have $i$ maximal elements,
\item[{\rm (ii)}] the number of bounded dominant regions $R$ of $\Cat_\Phi$
such that $i$ walls of $R$ of the form $H_{\alpha, 1}$ do not separate 
$R$ from $A_\circ$,
\item[{\rm (iii)}] the number of $w \in W (\Phi)$ such that $i$ walls of 
$w A_\circ$ of the form $H_{\alpha, 1}$ do not separate $w A_\circ$ from 
$A_\circ$ and
\item[{\rm (iv)}] the number of elements $w \in W (\Phi)$ having $i$ right 
ascents.
\end{enumerate}
\label{cor:shi}
\end{corollary}

\noindent
\emph{Proof.}
This follows from the case $m=1$ of Corollary \ref{cor:ik} and \cite[Lemma 
2.5]{Ath2}.
\qed

\section{Coroot lattice points and the affine Weyl group}
\label{proof}

In this section we complete the proof of Theorem \ref{thm1} (see Corollary 
\ref{cor:proof}). We assume that $\Phi$ is irreducible and crystallographic 
of rank $\ell$.

\begin{figure}
\epsfysize = 3.0 in 
\centerline{\epsffile{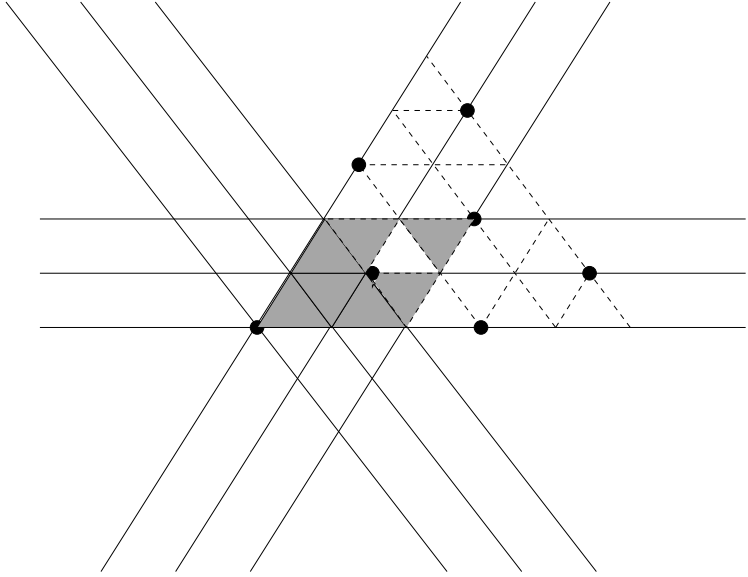}} 
\caption{The maximal alcoves of the bounded dominant regions and the 
simplex $p \, \overline{A_\circ}$ for $\Phi = A_2$ and $m=2$.}
\label{figure1}
\end{figure}

As in \cite[Section 4]{Ath2}, by the reflection in $W$ corresponding to 
a hyperplane $H_{\alpha, k}$ we mean the reflection in the linear hyperplane 
$H_\alpha$. We let $p = mh-1$, as in Section \ref{pre}, and $D_m (\Phi) = 
\check{Q} \cap p \, \overline{A_\circ}$. The following elementary lemma, 
for which a detailed proof can be found in \cite[Section 7.4]{Ha}, implies 
that $D_m (\Phi)$ is a set of representatives for the orbits of the 
$W$-action on $T_m (\Phi)$. 
\begin{lemma} \mbox{\rm (cf. \cite[Lemma 7.4.1]{Ha})}
The natural inclusion map from $D_m (\Phi)$ to the set $O_m (\Phi)$ of 
orbits of the $W$-action on $T_m (\Phi)$ is a bijection.

Moreover, if $y \in D_m (\Phi)$ then the stabilizer of $y$ with respect 
to the $W$-action on $T_m (\Phi)$ is the subgroup of $W$ generated by the 
reflections corresponding to the walls of $p \, \overline{A_\circ}$ 
which contain $y$. In particular, $r(y)$ is equal to the number of 
walls of $p \, \overline{A_\circ}$ which contain $y$. \qed
\label{lem:hai}
\end{lemma}
We will define a bijection $\rho: W_m (\Phi) \to D_m (\Phi)$ such 
that for $w \in W_m (\Phi)$, the number of walls of $w A_\circ$ of
the form $H_{\alpha, m}$ which do not separate $w A_\circ$ from $A_\circ$ 
is equal to the number of walls of $p \, \overline{A_\circ}$ which 
contain $\rho (w)$. Let $R_f$ be the region of $\aA^m_\Phi$ defined 
by the inequalities $m-1 < (\alpha, x) < m$ for $1 \le i \le \ell$. 
Let $w_f = w_{R_f}$ be the unique element $w$ of $W_m (\Phi)$ such 
that $w A_\circ \subseteq R_f$. We define the map $\rho: W_m (\Phi) 
\to \check{Q}$ by
\[ \rho (w) = (w_f \, w^{-1}) \cdot 0 \]
for $w \in W_m (\Phi)$. Observe that, by Lemma \ref{lem:h}, the alcove 
$w_f A_\circ$ can be described explicitly as the open simplex in $V$ 
defined by the linear inequalities $(\sigma_i, x ) < m$ for $1 \le i 
\le \ell$ and $(\tilde{\alpha}, x) > mh-m-1$. For any $1 \le r \le m$ 
we define the simplex
\[ \Sigma^r_m = \{x \in V: \, m-r \, \le \, (\sigma_i, x ) \ 
{\rm for} 
\ 1 \le i \le \ell \ {\rm and} \ (\tilde{\alpha}, x) \, \le \, 
mh-m+r-1\},\]
so that $\Sigma^m_m = p \, \overline{A_\circ}$. For any 
$\ell$-dimensional simplex $\Sigma$ in $V$ bounded by hyperplanes 
$H_{\alpha, k}$ in $\widetilde{\aA}_\Phi$ with $\alpha \in \Pi \cup 
\{\tilde{\alpha}\}$ we denote by $H(\Sigma, i)$ the wall of $\Sigma$ 
orthogonal to $\tilde{\alpha}$ or $\sigma_i$, if $i = 0$ or $i > 0$, 
respectively. We write $H(w, i)$ instead of $H(w \overline{A_\circ}, 
i)$ for $w \in W_a$. The reader is invited to test the results that 
follow in the case pictured in Figure \ref{figure1}. 
\begin{theorem}
The map $\rho$ is a bijection from $W_m (\Phi)$ to $D_m (\Phi)$. Moreover 
for any $w \in W_m (\Phi)$, $1 \le r \le m$ and $0 \le i \le \ell$, the 
point $\rho(w)$ lies on the wall $H(\Sigma^r_m, i)$ if and only if 
the wall $(w w_f^{-1}) \, H(w_f, i)$ of $w A_\circ$ is of the form 
$H_{\alpha, r}$ and does not separate $w A_\circ$ from $A_\circ$. 
\label{bij}
\end{theorem}

\noindent
\emph{Proof.}
Analogous to the proof of \cite[Theorem 4.2]{Ath2}.
\qed

\bigskip

\begin{corollary}
For any nonnegative integers $i_1, i_2,\dots,i_m$ each of the quantities 
which appear in the statement of Corollary \ref{cor:ik} is equal to the
number of points in $D_m (\Phi)$ which lie in $i_r$ walls of $\Sigma^r_m$ 
for all $1 \le r \le m$.
\label{ik:new}
\end{corollary}

\noindent
\emph{Proof.}
This follows from Theorem \ref{bij}.
\qed

\medskip
The next corollary completes the proof of Theorem \ref{thm1}.
\begin{corollary}
For any $0 \le i \le \ell$ the following are equal to $h^+_{\ell-i} 
(\Phi, m)$:
\begin{enumerate}
\itemsep=0pt
\item[{\rm (i)}] the number of points in $D_m (\Phi)$ which lie in $i$ walls 
of $p \, \overline{A_\circ}$ and  
\item[{\rm (ii)}] the number of $w \in W_m (\Phi)$ such that $i$ walls of 
$w A_\circ$ of the form $H_{\alpha, m}$ do not separate $w A_\circ$ 
from $A_{\circ}$.
\end{enumerate}
\label{cor:proof}
\end{corollary}

\noindent
\emph{Proof.}
By Lemma \ref{lem:hai}, the number of orbits $x \in O_m (\Phi)$ with 
rank $r(x) = i$ is equal to the number of points in $D_m (\Phi)$
which lie in $i$ walls of $p \, \overline{A_\circ}$. The statement now
follows by specializing Corollary \ref{ik:new} and recalling that
$\Sigma^m_m = p \, \overline{A_\circ}$.
\qed

\begin{remark} {\rm
Part (ii) of Corollary \ref{cor:proof} implies that $h^+_\ell (\Phi, m)$ 
is equal to the cardinality of $\check{Q} \cap \, (mh-1) A_\circ$. An argument 
similar to the one employed in \cite[Remark 4.5]{Ath2} shows that $\check{Q} 
\cap \, (mh-1) A_\circ$ is equinumerous to $\check{Q} \cap \, (mh-h-1) 
\overline{A_\circ}$. Therefore from Theorem \ref{thm0} (iii) we conclude 
that $h^+_\ell (\Phi, m) = N^+ (\Phi, m-1)$. Since the reduced Euler 
characteristic $\widetilde{\chi} (\Delta^m_+ (\Phi))$ of $\Delta^m_+ 
(\Phi)$ is equal to $(-1)^{\ell-1} h_\ell (\Delta^m_+ (\Phi))$, it follows 
from the results of the next section that $$\widetilde{\chi} (\Delta^m_+ 
(\Phi)) = (-1)^{\ell-1} N^+ (\Phi, m-1)$$ if $\Phi$ has type $A$, $B$ or 
$C$ (see \cite[(3.1)]{FR2} for the corresponding property of $\Delta^m 
(\Phi)$).  
\qed
\label{rem}} 
\end{remark}

The following conjecture is the positive analogue of \cite[Conjecture 
3.7]{FR2}.

\begin{conjecture} 
For any crystallographic root system $\Phi$ and all $m \ge 1$ the complex 
$\Delta^m_+ (\Phi)$ is pure $(\ell-1)$-dimensional and shellable and has 
reduced Euler characteristic equal to $(-1)^{\ell-1} N^+ (\Phi, m-1)$. 

In particular it is Cohen-Macaulay and has the homotopy type of a wedge of 
$N^+ (\Phi, m-1)$ spheres of dimension $\ell-1$.
\label{conj1}
\end{conjecture}

\section{The numbers $f_i (\Phi, m)$ and $f^+_i (\Phi, m)$}
\label{f}

Let $\Phi$ be a crystallographic root system of rank $\ell$ spanning the 
Euclidean space $V$. We define numbers $f_i (\Phi, m)$ and $f^+_i (\Phi, 
m)$ by the relations
\begin{equation}
\sum_{i=0}^\ell \, f_{i-1} (\Phi, m) (x-1)^{\ell-i} = \sum_{i=0}^\ell \, 
h_i (\Phi, m) \, x^{\ell-i} 
\label{f_i}
\end{equation}
and
\begin{equation}
\sum_{i=0}^\ell \, f^+_{i-1} (\Phi, m) (x-1)^{\ell-i} = \sum_{i=0}^\ell \, 
h^+_i (\Phi, m) \, x^{\ell-i} 
\label{f_i+}
\end{equation}
respectively. Comparing to equation (\ref{f+FR}) we see that Conjecture 
\ref{conj0} for the pair $(\Phi, m)$ is equivalent to the statement that
\begin{equation}
f^+_{i-1} (\Phi, m) = f_{i-1} (\Delta^m_+ (\Phi))
\label{conj00}
\end{equation}
for all $i$, where $f_{i-1} (\Delta^m_+ (\Phi))$ is as in Section \ref{pre}. 
We will give a combinatorial interpretation to the numbers $f _{i-1} (\Phi, 
m)$ and $f^+_{i-1} (\Phi, m)$ as follows. For $0 \le k \le \ell$ we denote 
by $\fF_k (\Phi, m)$ the collection of $k$-dimensional (nonempty) sets of 
the form 
\begin{equation}
\bigcap_{(\alpha, r) \, \in \, \Phi^+ \times \{0, 1,\dots,m\}} 
\ \tilde{H}_{\alpha, r}
\label{cells}
\end{equation}
where $\tilde{H}_{\alpha, r}$ can be 
\[ \cases{
H^+_{\alpha, 0}, & if \ $r=0$, \cr
H^-_{\alpha, m}, \, H^+_{\alpha, m} \ {\rm or} \ H_{\alpha, m}, & 
if \ $r=m$, \cr
H^-_{\alpha, r} \ {\rm or} \ H^+_{\alpha,r}, & if \ $1 \le r < m$ } \]
and $H^-_{\alpha, r}$ and $H^+_{\alpha, r}$ denote the two open half-spaces 
in $V$ defined by the inequalities $(\alpha, x) < r$ and $(\alpha, x) > r$, 
respectively. Observe that each element of $\fF_k (\Phi, m)$ is dominant. 
We also denote by $\fF^+_k (\Phi, m)$ the elements of $\fF_k (\Phi, m)$ 
which are bounded subsets of $V$ or, equivalently, the sets of the form 
(\ref{cells}) with $\tilde{H}_{\sigma_i, m} = H^-_{\sigma_i, m}$ or 
$H_{\sigma_i, m}$ for $1 \le i \le \ell$. In the special case $m=1$ part 
(iii) of the following theorem is the content of Remark 5.10 (v) in 
\cite{FR1}.
\begin{theorem}
For any irreducible crystallographic root system $\Phi$ and all $m \ge 1$ 
and $0 \le k \le \ell$ the number $f_{k-1} (\Phi, m)$ counts
\begin{enumerate}
\itemsep=0pt
\item[{\rm (i)}] pairs $(R, S)$ where $R$ is a dominant region of $\aA_\Phi^m$ and 
$S$ is a set of $\ell-k$ walls of $R$ of the form $H_{\alpha, m}$ which 
separate $R$ from $A_\circ$,
\item[{\rm (ii)}] pairs $(\iI, T)$ where $\iI$ is a geometric chain of filters in 
$\Phi^+$ of length $m$ and $T$ is a set of $\ell-k$ indecomposable roots of 
rank $m$ with respect to $\iI$ and 
\item[{\rm (iii)}] the elements of $\fF_k (\Phi, m)$.
\end{enumerate}
\label{thm:eleni1}
\end{theorem}

\noindent
\emph{Proof.}
>From (\ref{f_i}) we have
\[ f_{k-1} (\Phi, m) = \sum_{i=0}^k \ h_i (\Phi, m) {\ell-i \choose \ell-k} 
\]
which cleary implies (i) and (ii) (see Theorem \ref{thm:ath2}). To complete 
the proof it suffices to give 
a bijection from the set $\rR_k (\Phi, m)$ of pairs $(R, S)$ which appear in 
(i) to $\fF_k (\Phi, m)$. Given such a pair $\tau = (R, S)$ let $g(\tau)$ be 
the intersection (\ref{cells}), where $\tilde{H}_{\alpha, r}$ is chosen so 
that $R \subseteq \tilde{H}_{\alpha, r}$ unless $r=m$ and $H_{\alpha, r} \in 
S$, in which case $\tilde{H}_{\alpha, r} = H_{\alpha, r}$. Let $S = 
\{H_{\alpha_1, m}, H_{\alpha_2, m},\dots,H_{\alpha_{\ell-k}, m}\}$ and let
$F_S$ be the intersection of the hyperplanes in $S$. It follows from 
\cite[Corollary 3.14]{Ath2} that $S$ is a proper subset of the set of walls
of an alcove of $\widetilde{\aA}_\Phi$ and hence that $F_S$ is nonempty and 
$k$-dimensional. To show that $g(\tau)$ is nonempty and $k$-dimensional, so 
that $g: \rR_k (\Phi, m) \to \fF_k (\Phi, m)$ is well defined, we need to 
show that $F_S$ is not contained in any hyperplane $H_{\alpha, r}$ with 
$\alpha \in \Phi^+$ and $0 \le r \le m$ other than those in $S$. So suppose
that $F_S \subseteq H_{\alpha, r}$ with $\alpha \in \Phi^+$ and $r \ge 0$.
Then there are real numbers $\lambda_1, \lambda_2,\dots,\lambda_{\ell-k}$
such that 
\begin{equation}
\alpha = \lambda_1 \alpha_1 + \lambda_2 \alpha_2 + \cdots + \lambda_{\ell-k} 
\alpha_{\ell-k}
\label{lambdas}
\end{equation}
and $r = m (\lambda_1 + \lambda_2 + \cdots + \lambda_{\ell-k})$. Observe 
that the $\alpha_i$ are minimal elements of the last filter in the geometric 
chain of filters in $\Phi^+$ corresponding to $R$ and hence that they form 
an antichain in $\Phi^+$, meaning a set of pairwise incomparable elements. 
It follows from the first main result of \cite{So} (see the proof of 
\cite[Corollary 6.2]{AR}) that the coefficients $\lambda_i$ in 
(\ref{lambdas}) are nonnegative integers. Hence either $r>m$ or $\alpha 
= \alpha_i$ and $r=m$ for some $i$, so that $H_{\alpha, r} \in S$.

To show that $g$ is a bijection we will show that given $F \in \fF_k (\Phi, 
m)$ there exists a unique $\tau \in \rR_k (\Phi, m)$ with $g(\tau)=F$. Let
$(\varpi_1^\vee, \varpi_2^\vee,\dots,\varpi_\ell^\vee)$ be the linear 
basis of $V$ which is dual to $\Pi$, in the sense that
\[ (\sigma_i, \varpi_j^\vee) = \delta_{ij}. \]
Observe that if $g(\tau) = F$ with $\tau = (R, S)$, $x$ is a point in $F$ 
and $\epsilon_i$ are sufficiently small positive numbers then
\begin{equation}
x + \sum_{i=1}^\ell \ \epsilon_i \varpi_i^\vee \, \in \, R.
\label{epsilon}
\end{equation}
Since regions of $\aA^m_\Phi$ are pairwise disjoint this implies uniqueness 
of $R$, and hence of $\tau$. To prove the existence let $R$ be the unique
region of $\aA^m_\Phi$ defined by (\ref{epsilon}). Equivalently $R$ can be
obtained by replacing all hyperplanes of the form $H_{\alpha, m}$ in the 
intersection (\ref{cells}) defining $F$ by $H^+_{\alpha, m}$. It suffices to 
show that any such hyperplane $H_{\alpha, m}$ is a wall of $R$ since then,
if $S$ is the set of hyperplanes of the form $H_{\alpha, m}$ which contain
$F$ then $\tau = (R, S) \in \rR_k (\Phi, m)$ and $g (\tau) = F$. 

Suppose on the contrary that $H_{\alpha, m} \supseteq
F$ is not a wall of $R$. It follows from Theorem \ref{thm:ath2} that 
$\alpha$ is not indecomposable of rank $m$ with respect to the geometric 
chain $\iI$ of filters $\Phi^+ = \iI_0 \supseteq \iI_1 \supseteq \cdots 
\supseteq \iI_m$ corresponding to $R$ and hence that one can write $\alpha =
\beta + \gamma$ for some $\beta \in \iI_i$, $\gamma \in \iI_j$ with $i+j = m$.
Since $(\alpha, x) = m$ for $x \in F$ and $(\beta, x) > i$ and $(\gamma, x) 
> j$ hold for $x \in R$, so that $(\beta, x) \ge i$ and $(\gamma, x) \ge j$ 
hold for $x \in F$, we must have $(\beta, x) = i$ and $(\gamma, x) = j$ 
for $x \in F$. However one of $i,j$ must be less than $m$ and this 
contradicts the fact that $F$ can be contained in $H_{\alpha, r}$ for 
$\alpha \in \Phi^+$ only if $r=m$.
\qed

\medskip
The proof of the next theorem is entirely similar to that of Theorem 
\ref{thm:eleni1} and is omitted.
\begin{theorem}
For any irreducible crystallographic root system $\Phi$ and all $m \ge 1$ 
and $0 \le k \le \ell$ the number $f^+_{k-1} (\Phi, m)$ counts
\begin{enumerate}
\itemsep=0pt
\item[{\rm (i)}] pairs $(R, S)$ where $R$ is a dominant bounded region of 
$\aA_\Phi^m$ and $S$ is a set of $\ell-k$ walls of $R$ of the form $H_{\alpha, 
m}$ which do not separate $R$ from $A_\circ$,
\item[{\rm (ii)}] pairs $(\jJ, T)$ where $\jJ$ is a positive geometric chain of 
ideals in $\Phi^+$ of length $m$ and $T$ is a set of $\ell-k$ indecomposable 
roots of rank $m$ with respect to $\jJ$ and
\item[{\rm (iii)}] the elements of $\fF^+_k (\Phi, m)$. \qed
\end{enumerate}
\label{thm:eleni2}
\end{theorem}

The reader is invited to use part (iii) of Theorems \ref{thm:eleni1} and 
\ref{thm:eleni2} as well
as Figure \ref{figure1} to verify that $f_{-1} = 1$, $f_0 = 8$, $f_1 = 12$,
$f^+_{-1} = 1$, $f^+_0 = 6$ and $f^+_1 = 7$ when $\Phi = A_2$ and $m=2$.
It should be clear that apart from the necessary modifications in the 
statements of part (i), Theorems \ref{thm:eleni1} and \ref{thm:eleni2} 
are also valid when $\Phi$ is reducible.
\begin{lemma}
For any crystallographic root system $\Phi_I$ and all $m \ge 1$ and $0 \le 
k \le \ell$ we have
\[ f_{k-1} (\Phi_I, m) = \sum_{J \subseteq I} \ f^+_{k-|J|-1} (\Phi_{I 
\sm J}, m). \]
\label{lem:ff+}
\end{lemma}

\noindent
\emph{Proof.}
For $J \subseteq I$ let $V_J$ be the linear span of the simple roots indexed 
by the elements of $J$ and let $p_J: V_I \to V_{I \sm J}$ be the orthogonal 
projection onto $V_{I \sm J}$. We define the \emph{simple part} of $F \in 
\fF_k (\Phi_I, m)$ as the set of indices $j \in I$ such that $F \subseteq 
H^+_{\sigma_j, m}$. Observe that if $J$ is the simple part of $F$ then for 
$\alpha \in \Phi^+$ and $x \in F$ we have 
\[ \begin{tabular}{ll} $(\alpha, x) = (\alpha, p_J (x))$, & 
{\rm if} \ $\alpha \in \Phi_{I \sm J}$ 
\\ $(\alpha, x) > m$, & {\rm otherwise.}
\end{tabular} \]
It follows that $p_J$ induces a bijection from the set of elements of $\fF_k 
(\Phi_I, m)$ with simple part $J$ to $\fF^+_{k-|J|} (\Phi_{I \sm J}, m)$. 
Hence counting the elements of $\fF_k (\Phi_I, m)$ according to their simple 
part proves the lemma.
\qed

\medskip
The same type of argument as that in the next corollary appears in the proof
of \cite[Proposition 6.1]{MRZ}.

\begin{corollary} 
If for some pair $(\Phi, m)$ every parabolic subsystem of $\Phi$ satisfies 
{\rm (\ref{ff+FR})} and we have $h_i (\Phi, m) = h_i (\Delta^m (\Phi))$ for 
all $i$ then we have $h^+_i (\Phi, m) = h_i (\Delta^m_+ (\Phi))$ for all $i$ 
as well.
\label{cor:if}
\end{corollary}

\noindent
\emph{Proof.}
>From the assumption we have $f_{k-1} (\Phi, m) = f_{k-1} (\Delta^m (\Phi))$ 
for all $k$. Equation (\ref{ff+FR}) and Lemma \ref{lem:ff+} imply that 
$f^+_{k-1} (\Phi, m) = f_{k-1} (\Delta^m_+ (\Phi))$ for all $k$ via M\"obius 
inversion on 
the set of pairs $(k, I)$ partially ordered by letting $$(l, J) \le (k, I) 
\ {\rm if \ and \ only \ if} \ J \subseteq I \ {\rm and} \ k-l = |I \sm J|.$$  
This is equivalent to the conclusion of the corollary.
\qed

\begin{corollary} 
Conjecture \ref{conj0} holds for root systems of type $A$, $B$ or $C$ and any 
$m\ge 1$ and for all root systems when $m=1$.
\label{cor:conj}
\end{corollary}

\noindent
\emph{Proof.}
This follows from Lemma \ref{lem:FR}, the previous corollary and the fact 
that the equality $h_i (\Phi, m) = h_i (\Delta^m (\Phi))$ can be checked 
case by case from the explicit formulas given in \cite{Ath2, FR2, Tz} in 
the cases under consideration.
\qed

\medskip
We conclude this section with a combinatorial interpretation to $f^+_{k-1} 
(\Phi, m)$ similar to those provided in parts (i) and (ii) of Theorem 
\ref{thm:eleni2}. 
\begin{theorem}
For any irreducible crystallographic root system $\Phi$ and all $m \ge 1$ 
and $0 \le k \le \ell$ the number $f^+_{k-1} (\Phi, m)$ counts
\begin{enumerate}
\itemsep=0pt
\item[{\rm (i)}] pairs $(R, S)$ where $R$ is a dominant region of $\aA_\Phi^m$ 
and $S$ is a set of $\ell-k$ walls of $R$ of the form $H_{\alpha, 
m}$ which separate $R$ from $A_\circ$ such that $S$ contains all such 
walls of $R$ with $\alpha \in \Pi$ and
\item[{\rm (ii)}] pairs $(\iI, T)$ where $\iI$ is a geometric chain of filters in 
$\Phi^+$ of length $m$ and $T$ is a set of $\ell-k$ indecomposable roots of 
rank $m$ with respect to $\iI$ which contains all simple indecomposable 
roots of rank $m$ with respect to $\iI$.
\end{enumerate}
\label{thm:eleni3}
\end{theorem}

\noindent
\emph{Proof.}
The sets in (i) and (ii) are equinumerous by Theorem \ref{thm:ath2}.
To complete the proof one can argue that the map $g$ in the proof of 
Theorem \ref{thm:eleni1} restricts to a bijection from the set in (i) to 
$\fF^+_k (\Phi, m)$. Alternatively, arguing as in the proof of 
Corollary \ref{cor:if}, it suffices to show that
\[ f_{k-1} (\Phi_I, m) = \sum_{J \subseteq I} \ 
g^+_{k-|J|-1} (\Phi_{I \sm J}, m), \]
where $\Phi = \Phi_I$ and $g^+_{k-1} (\Phi, m)$ denotes the cardinality 
of the set of pairs, say $\gG^+_k (\Phi, m)$, which appears in (ii). Let 
$\gG_k (\Phi_I, m)$ denote the set of pairs defined in (ii) of the 
statement of Theorem 
\ref{thm:eleni1}. For $(\iI, T) \in \gG_k (\Phi_I, m)$ call the set of 
simple roots which are indecomposable of rank $m$ with respect to $\iI$ 
and are not contained in $T$ the \emph{simple part} of $(\iI, T)$ and 
for any $J \subseteq I$ denote by $\Lambda_J$ the order filter of roots 
$\alpha \in \Phi^+_I$ for which $\sigma \le \alpha$ for some $\sigma \in 
J$. It is straightforward to check from the definitions that the map 
which sends a pair $(\iI, T)$ to $(\iI \sm \Lambda_J, T)$, where $\iI 
\sm \Lambda_J$ denotes the chain obtained from $\iI$ by removing 
$\Lambda_J$ from each filter of $\iI$, induces a bijection from the set 
of elements of $\gG_k (\Phi_I, m)$ with simple part $J$ to $\gG^+_{k-|J|} 
(\Phi_{I \sm J}, m)$. Therefore counting the elements of $\gG_k (\Phi_I, 
m)$ according to their simple part gives the desired equality.  
\qed

\section{Classical types and the case $m=1$}
\label{class}

In this section we compute the numbers $h^+_i (\Phi, m)$ and $f^+_{i-1} 
(\Phi, m)$ in the cases of the classical root systems. 

\begin{proposition} 
The number $h^+_i (\Phi, m)$ is equal to 

\begin{center}
\begin{tabular}{lr}
{\Large $\frac{1}{i+1} {n-1 \choose i} {mn-2 \choose i}$}, & {\rm if} \ 
$\Phi = A_{n-1}$, \\ \\
{\Large ${n \choose i} {mn-1 \choose i}$}, & {\rm if} \ $\Phi = B_n$ 
{\rm or} $C_n$, \\ \\
{\Large ${n \choose i} {m(n-1)-1 \choose i} + {n-2 \choose i-2} 
{m(n-1) \choose i}$}, & {\rm if} \ $\Phi = D_n$.
\end{tabular}
\end{center}

\label{thm:ABCD}
\end{proposition}

\noindent
\emph{Proof.}
The proof can be obtained from that of \cite[Proposition 5.1]{Ath2} by 
replacing the quantity $mh+1$ by $mh-1$ and using Theorem 1.3 (iii) 
instead of \cite[Theorem 1.2 (ii)]{Ath2}.  
\qed

\medskip
The following corollary is a straightforward consequence of Proposition 
\ref{thm:ABCD} and equation (\ref{f_i+}).

\begin{corollary} 
The number $f^+_{k-1} (\Phi, m)$ is equal to  

\begin{center}
\begin{tabular}{lr}
{\Large $\frac{1}{k+1} {n-1 \choose k} {mn+k-1 \choose k}$}, & {\rm if} \ 
$\Phi = A_{n-1}$, \\ \\
{\Large ${n \choose k} {mn+k-1 \choose k}$}, & {\rm if} \ $\Phi = 
B_n$ {\rm or} $C_n$, 
\\ \\
{\Large ${n \choose k} {m(n-1)+k-1 \choose k} + {n-2 \choose k-2} 
{m(n-1)+k-2 \choose k}$}, & {\rm if} \ $\Phi = D_n$.
\end{tabular}
\end{center}
\qed
\label{cor:ABCD}
\end{corollary}

In the case $m=1$, the following corollary for $\Phi = A_{n-1}$ and 
$\Phi = B_n, C_n$ is a special case of \cite[(34)]{Ch} and \cite[(46)]{Ch}, 
respectively.
\begin{corollary} 
The number $f_{k-1} (\Delta^m_+(\Phi))$ is equal to  

\begin{center}
\begin{tabular}{lr}
{\Large $\frac{1}{k+1} {n-1 \choose k} {mn+k-1 \choose k}$}, & {\rm if} \ 
$\Phi = A_{n-1}$, \\ \\
{\Large ${n \choose k} {mn+k-1 \choose k}$}, & {\rm if} \ $\Phi = B_n$ 
{\rm or} $C_n$.
\end{tabular}
\end{center}
Moreover $$f_{k-1} (\Delta_+(D_n)) = {n \choose k} {n+k-2 \choose k} + 
{n-2 \choose k-2} {n+k-3 \choose k}.$$
\label{cor:AB}
\end{corollary}

\noindent
\emph{Proof.}
Combine Corollaries \ref{cor:conj} and \ref{cor:ABCD}.
\qed

\begin{remark} {\rm
The number of positive filters in $\Phi^+$ with $i$
minimal elements has been computed for the exceptional root systems
by Victor Reiner as shown in the following table.} 
\end{remark}

{\scriptsize
\begin{table}[hptb]
\begin{center}
\begin{tabular}{| l| l | l | l | l | l | l | l | l | l |} \hline
 \ \ \ \ $i$ & 0 & 1  & 2 & 3 & 4 & 5 & 6 & 7 \\ \hline \hline
 $\Phi = G_2$  &1    & 4  &  &  &  &  &  &    \\ \hline
 $\Phi = F_4$  &1    & 20 &  35 &  10 &  &  &  &  \\ \hline
 $\Phi = E_6$  &1    & 30 & 135 & 175 & 70 & 7 &  &  \\ \hline
 $\Phi = E_7$  &1    & 56 &  420 & 952 & 770 & 216 & 16 & \\ \hline
 $\Phi = E_8$  &1    & 112 & 1323 & 4774 & 6622 & 3696 & 770 & 44  
\\ \hline
\end{tabular}
\caption{The numbers $h^+_i (\Phi)$ for the exceptional root systems.}
 \label{except}
\end{center}
\end{table}
}

\vspace{0.1 in}
\noindent
\emph{Proof of Theorem \ref{thm2}.}
We will prove the statement of the theorem without the assumption that 
$\Phi$ is irreducible. Let $\ell$ be the rank of $\Phi = \Phi_I$. We 
write $h_k (\Phi_I)$ instead of $h_k (\Phi_I, 1)$, so that $h_{\ell-k} 
(\Phi_I)$ counts the filters in $\Phi_I^+$ with $k$ minimal elements as 
well as the ideals in $\Phi_I^+$ with $k$ maximal elements. Let 
$\widetilde{h}^+_k (\Phi)$ denote the number of positive filters in 
$\Phi^+$ with $k$ minimal elements. Counting filters in $\Phi_I^+$ by 
the set of simple roots they contain gives 
\begin{equation}
h_{\ell-k} (\Phi_I) = \sum_{J \subseteq I} \ \widetilde{h}^+_{k-|J|} 
(\Phi_{I \sm J}). 
\label{dual1}
\end{equation}
Similarly, counting ideals in $\Phi_I^+$ by the set of simple roots 
they do not contain gives 
\[ h_{\ell-k} (\Phi_I) = \sum_{J \subseteq I} \ h^+_{\ell-|J|-k} 
(\Phi_{I \sm J}). \]
Since it is known \cite{Pa} that $h_{\ell - k} (\Phi_I) = h_k (\Phi_I)$
the previous relation can also be written as
\begin{equation}
h_{\ell-k} (\Phi_I) = \sum_{J \subseteq I} \ h^+_{k - |J|} (\Phi_{I \sm 
J}). 
\label{dual2}
\end{equation}
Comparing (\ref{dual1}) and (\ref{dual2}) and using M\"obius inversion 
as in Section \ref{f} gives $\widetilde{h}^+_k (\Phi_J) = h^+_k (\Phi_J)$ 
for all $J \subseteq I$, which is the statement of the theorem. 
\qed

\section{Remarks}
\label{remarks}

\noindent
1. The following reciprocity relation
\begin{equation}
N^+ (\Phi, m-1) = (-1)^\ell N (\Phi, -m)
\label{rec}
\end{equation}
was observed by Fomin and Reading \cite[(2.12)]{FR2}. We will show that, 
as suggested by S. Fomin (private communication), this relation is in 
fact an instance of Ehrhart reciprocity. Let $i(n)$ be the cardinality
of $\check{Q} \cap \, n \overline{A_\circ}$ for $n \in \NN$. It is clear 
that the vertices of the simplex $\overline{A_\circ}$ have rational 
coordinates in the basis $\Pi$ of $V$ of simple roots, hence also in the 
basis $\Pi^\vee = \left\{ 2 \alpha / (\alpha,\alpha): \ \alpha \in 
\Pi \right\}$ of $V$. Therefore the function $i(n)$ is the Ehrhart 
quasi-polynomial of $\overline{A_\circ}$ with respect to the lattice 
$\check{Q}$ (see \cite[Section 4.6]{Sta} for an introduction to the 
theory of Ehrhart quasi-polynomials). Ehrhart reciprocity \cite[Theorem 
4.6.26]{Sta} implies that
\[ (-1)^\ell i(-n) = \# \, (\check{Q} \cap \, n A_\circ) \]
for $n \in \NN$ and hence, setting $n=mh-1$ and consulting \cite[Theorem 
1.1]{Ath2}, that
\[ (-1)^\ell N (\Phi, -m) = \# \, (\check{Q} \cap \, (mh-1) A_\circ). \]
Remark \ref{rem} asserts that
\[ \# \, (\check{Q} \cap \, (mh-1) A_\circ) = N^+ (\Phi, m-1) \]
and hence (\ref{rec}) holds.

\vspace{0.1 in}
\noindent
2. It would be interesting to give combinatorial proofs of the formulas 
in Corollary \ref{cor:AB} directly from the description of the relevant
complexes given in \cite{FR2, FZ, Tz} and Section \ref{pre}.

\vspace{0.1 in}
\noindent
3. After this paper was completed the following came to our attention. 
(i) The numbers $\widetilde{h}^+_i (\Phi)$ of positive filters in $\Phi^+$ 
with $i$ minimal elements are also discussed and partially computed in 
\cite[Section 3]{Pa2}. (ii) Theorem 2.7 in \cite{FR3} implies that 
Lemma \ref{lem:FR} is valid for all pairs $(\Phi, m)$. In view of (ii) 
and the equality $h_i (\Phi, m) = h_i (\Delta^m (\Phi))$ (see \cite{Ath2, 
FR2, FR3}) when $\Phi = D_n$, it follows from Corollary 
\ref{cor:if} that Conjecture \ref{conj0} is also valid for root systems 
of type $D$ and arbitrary $m$.  

\vspace*{0.25 in}
\noindent
\emph{Acknowledgements}. We are grateful to Sergey Fomin and Nathan 
Reading for making their work \cite{FR2} available to us. We also thank 
Sergey Fomin for useful discussions and Victor Reiner for providing 
the data of Table \ref{except} as well as the software to confirm it.

\end{document}